\documentclass[11pt]{article}
\usepackage[margin=1in]{geometry}
\usepackage{times}
\usepackage{amsmath,amsthm,amssymb,amsfonts,bm}
\usepackage{color}
\usepackage[colorlinks,linkcolor=black,bookmarksopen,
bookmarksnumbered,citecolor=black,urlcolor=black]{hyperref}
\usepackage{graphicx}

\newcommand{\nchoosek}[2]{\left(\begin{array}{c}#1\\#2\end{array}\right)}
\newcommand{\vd}{ \mathbf{d} }
\newcommand{\vep}{\bm{\epsilon}}
\newcommand{\vw}{ \mathbf{w} }
\newcommand{\vW}{ \mathbf{W} }
\newcommand{\vone}{ \mathbf{1} }

\newcommand{\Seq}{\mathcal S}

\makeatletter
\def\imod#1{\allowbreak\mkern2mu({\operator@font mod}\,\,#1)}
\makeatother

\definecolor{darkgrn}{rgb}{0, 0.8, 0}

\theoremstyle{definition}

\theoremstyle{remark}

\begin{document}

\title{A Knapsack-like Code Using Recurrence Sequence Representations}

\author{\large Nathan Hamlin
\quad\quad
Bala Krishnamoorthy
\quad\quad William Webb \\
\normalsize
Department of Mathematics, Washington State University, Pullman,
    WA, USA.  \\
{\tt \{nhamlin,bkrishna,webb\}@math.wsu.edu}
}

\date{}

\maketitle

\begin{abstract}
\noindent We had recently shown that every positive integer can be
represented uniquely using a recurrence sequence, when certain
restrictions on the digit strings are satisfied. We present the
details of how such representations can be used to build a
knapsack-like public key cryptosystem. We also present new disguising
methods, and provide arguments for the security of the code against
known methods of attack.
\end{abstract}

\section{Introduction} \label{sec-intro}

One of the first public key cryptosystems was the traditional knapsack
code proposed by Merkle and Hellman \cite{MH78}. This code has the
advantage of fast encoding and decoding. Also, in more recent
developments, it has been shown that a quantum computer could make
factoring large numbers fast enough to make the RSA code insecure
\cite{Sh1997,Luetal2012}. However, it appears that quantum computers
would still struggle to make knapsack problems fast to solve
\cite{ArSc2003}.

Unfortunately, the traditional knapsack code was broken by two
different approaches -- by reversing the disguising steps \cite{S84}
and by a direct attack using lattice-based approaches \cite{LO85}. We
describe a new type of knapsack-like code along with new disguising
techniques, which make the code resistant to both these classes of
attacks. We first give a brief description of the traditional knapsack
code and its weaknesses.

\subsection{The traditional knapsack code} \label{ssec-tradknap}

The plaintext message is assumed to be an integer $M$, $0 \leq M <
2^n$. We consider the representation of $M$ in base $2$:
\begin{equation} \label{eq-Minbase2}
M = \sum_{i=0}^{n-1} \epsilon_i 2^i, ~~ 0 \leq \epsilon_i \leq 1.
\end{equation}
The creator of the code chooses a secret, superincreasing sequence
$\{s_i\}$, i.e., where $s_i > s_{i-1} + \cdots + s_0$. The secret
$s_i$ are then disguised by one or more modular multiplications of the
form
\begin{equation} \label{eq-modmult}
w_i = k s_i \imod{m},
\end{equation}
where $k$ and $m$ are kept secret. The $w_i$ are made public. The
sender of the message $M$ computes
\begin{equation} \label{eq-encodedM}
T = \sum_{i=0}^{n-1} \epsilon_i w_i.
\end{equation}
The encoded message $T$ is then sent over a possibly insecure
channel. The hope was that only the creator of the code, who knows $k,
m,$ and the $s_i$, can solve Equation (\ref{eq-encodedM}) for the
coefficients $\epsilon_i$. In particular, the disguising step given in
Equation (\ref{eq-modmult}) is easily reversed.

Shamir was able to break this code by calculating $k$ and $m$ and thus
reversing the disguising \cite{S84}. Although the $w_i$ appear random,
the fact that the $s_i$ are superincreasing can be exploited to yield
enough clues to determine $k$ and $m$, or at least equally useful
alternative values.

The other attack, which has received the most attention, tries to
solve Equation (\ref{eq-encodedM}) directly using lattice-based
approaches \cite{LO85,CJLOSS92}. In essence, the solutions to Equation
(\ref{eq-encodedM}) correspond to a lattice of vectors. Basis
reduction algorithms are efficient methods to find short vectors in,
and short bases of, lattices \cite{LLL82,LLS90}. Although not
guaranteed to do so, these approaches often find a shortest vector in
the given lattice \cite{GamNgupredlatred08}. A weakness in using the
base-$2$ representation of the message $M$ is that the corresponding
vector $\vep = \begin{bmatrix} \epsilon_0 & \epsilon_1 & \cdots &
  \epsilon_{n-1} \end{bmatrix}^T$ from Equations (\ref{eq-Minbase2})
and (\ref{eq-encodedM}) is likely to be the unique shortest vector in
an appropriately defined basis.

For example, if $n=1000$, the size of $M$ is roughly $2^{1000}$. After
using the $s_i$ and disguising to create the $w_i$ (using Equation
(\ref{eq-modmult})), the target sum $T$ might be of size $2^{1030}$
(depending on the number of disguising steps used). The expected
length of the true decoding vector $\vep = \begin{bmatrix} \epsilon_0
  & \epsilon_1 & \cdots & \epsilon_{999} \end{bmatrix}^T$ is $500$
(for simplicity, we talk about the {\em square} of the length
here). The number of $0$--$1$ $1000$-vectors of length less than $500$
is roughly $2^{999}$. Hence the probability of any one of them
equaling $T$ is $2^{-31}$. In practice, approaches using basis
reduction may find vectors with some $\epsilon_i=-1$ or
$\epsilon_i=2$, and so on. We may penalize such solutions in these
approaches, but cannot forbid them. But allowing such $\epsilon_i$
values tend to produce longer vectors. Thus it is indeed likely that
the shortest vector in the lattice is the true decoding vector. In our
new knapsack-like code, we try to prevent such lattice-based
approaches from being effective.

\section{A Code Using Recurrence Sequence Representations} \label{sec-coderecseq}

Hamlin and Webb recently presented a description of how to find a
unique representation of any positive integer using a recurrence
sequence $\{u_i\}$ as a base \cite{HaWe2012}. We now show how to
create a knapsack-like public key code using such representations. We
illustrate our construction on a small example, and refer the reader
to the above paper \cite{HaWe2012} for the complete proofs of all
assertions.

\subsection{An example code} \label{ssec-exmplcode}

Let $u_i$ satisfy the recurrence $u_{i+5} = u_{i+4} + u_{i+2} +
2u_{i+1} + 7u_i$ (the initial values are not vital, and we could take
them to be the standard ones, i.e., $u_0=1, u_1=1, u_2=1, u_3=2$, and
$u_4=4$). The signature is $\Seq = 10127$. The representation of any
positive integer $M$ is of the form
\begin{equation} \label{eq-Mrepr}
M = \sum_{i=0}^{n-1} d_i u_i,
\end{equation}
where the string of digits $d_{n-1} d_{n-2} \dots d_1 d_0$ must be
composed of blocks of digits which are lexicographically smaller than
$\Seq$. In this case, the allowed blocks of digits are $0, 100, 1010,
1011$, $10120, 10121$, $10122, 10123, 10124, 10125, 10126$. Hence, for
instance, $~1011|1023|0|1010~$ is a legitimate string, but
$~1010|110123100~$ is not. Notice that no allowed block begins with
$11$.

We now illustrate how to calculate this type of representation of any
number $M$ using a greedy approach. Although this calculation is
straightforward, the fact that makes this code harder to model for the
cryptanalyst is that so many strings of digits are not allowed in the
representations, even though the strings appear similar to the allowed
ones, and both classes of strings use digits of the same size.

\subsubsection{Computing the representation of $M$} \label{sssec-comprepM}
An easy way to calculate the representation of any number $M$ in the
recurrence sequence $\{u_i\}$ is to calculate the augmented sequence
$\{v_{j,i}\}$ for $1 \leq j \leq 10$, $0 \leq i \leq n$ given by
\begin{equation*} \label{eq-vij}
\begin{array}{rcl}
v_{1,i} & = & u_i, \\
v_{2,i} & = & u_i + u_{i-2}, \\
v_{3,i} & = & u_i + u_{i-2} + u_{i-3}, \\
v_{4,i} & = & u_i + u_{i-2} + 2u_{i-3}, \\
v_{5,i} & = & u_i + u_{i-2} + 2u_{i-3} + u_{i-4}, \\
& \vdots & \\
v_{10,i} & = & u_i + u_{i-2} + 2u_{i-3} + 6u_{i-4}.
\end{array}
\end{equation*}
The $v_{j,i}$ correspond to the allowed blocks of digits. In our
example, the $v_{j,i}$ occur in groups of size $10$. In other
examples, the groups could be much larger.

The correct expression for $M$ is found simply by using the greedy
algorithm on the $v_{j,i}$, and converting the sum into an expression
in the $u_i$. The $v_{j,i}$ could be calculated and stored, or
calculated from the $u_i$ as needed.

We explore the memory requirements for storing all the $\{v_{j,i}\}$.
The principal eigenvalue of the sequence $\{u_i\}$ is $\alpha \approx
1.9754$. We may assume that $u_i$ is roughly $\alpha^i$, or is close
to $2^i$. After disguising, the public weights $w_i$ will be
approximately $2^{n+40}$, and the target $T$ approximately
$2^{n+50}$. In other words, these quantities require $40$--$50$ extra
bits of memory to represent. If $n=1000$, the memory required for the
weights $w_i$ is roughly $1,040,000$ bits (or $130$ kilobytes). The
memory required to store all the $v_{j,i}$ is hence $1.3$
megabytes. Even if $n$ is much larger, the memory needed for the $w_i$
is negligible.

\subsection{Encryption and decryption} \label{ssec-encrdecr}

Let $\{u_i\}$ be a recurrence sequence that satisfies the following
recurrence equation.
\begin{equation} \label{eq-unreqeq}
u_i = a_1 u_{i-1} + a_2 u_{i-2} + \cdots + a_h u_{i-h},
\end{equation}
where $a_1 > 0$ and all $a_i \geq 0$. The string $\Seq = a_1a_2 \cdots
a_h$ is its signature, and we let $A = a_1 + a_2 + \cdots +
a_h$. Every natural number $N$ has a unique representation in the form
of Equation (\ref{eq-Mrepr}), where the digits are composed of blocks
that are lexicographically smaller than $\Seq$. Including the zero
block, there are $A$ such blocks. The auxiliary sequence $\{v_{j,i}\}$
is constructed as linear combinations of the $u_i$ with coefficients
same as the blocks other than the zero block. Hence there are $A-1$ of
the $v_{j,i}$ in each group. The total number of $v_{j,i}$ numbers is
hence $(A-1)n$ if there are $n$ of the $u_i$.

The creator of the code chooses a secret sequence $\{s_i\}$ which has
the property that $s_{i+1} \, > s_i \,(u_{i+1}/u_i)$ for all $i$. This
property replaces the condition of $\{s_i\}$ being superincreasing as
used in the traditional knapsack code.

The $s_i$ are then disguised by any invertible mapping, some of which
we describe below. The resulting quantities $w_i$ are the public
weights. If $M$ is the original plaintext message, the user of the
code expresses
\begin{equation} \label{eq-Mrepruivij}
M = \sum_{i=0}^{n-1} d_i u_i = \sum_{i=0}^{n-1} \epsilon_{j,i} \,
v_{j,i},
\end{equation}
and computes
\begin{equation} \label{eq-T}
T = \sum_{i=0}^{n-1} d_i w_i = \sum_{i=0}^{n-1} \epsilon_{j,i} \,
y_{j,i},
\end{equation}
which is the transmitted message.

\medskip
Since the mappings used for disguising are invertible by the code's
creator, she can calculate
\begin{equation} \label{eq-N}
N = \sum_{i=0}^{n-1} d_i s_i = \sum_{i=0}^{n-1} \epsilon_{j,i} \,
t_{j,i}.
\end{equation}
We must show that she can solve Equation (\ref{eq-N}) for the {\em
  same} digits $d_i$ as appearing in Equation
(\ref{eq-Mrepruivij}). The $t_{j,i}$ and $y_{j,i}$ are combinations of
the $s_i$ and $w_i$, respectively, in the same way as $v_{j,i}$ are
combinations of the $u_i$. Also, each $\epsilon_{j,i} = 0$ or $1$, and
for each $i$, at most one $\epsilon_{j,i} = 1$.

In the greedy algorithm to express $N$ using the $t_{j,i}$, we
subtract the largest possible $t_{j,i}$ from $N$ at each step, and
repeat the process on the remainder. If the correct $t_{j,i}$ have
been used previously, we have an equation at each step that is
essentially of the same form as Equation (\ref{eq-N}). That is, we
know the number $N'$ and that
\begin{equation} \label{eq-Npr}
N' = \sum_{i=0}^{k}  \epsilon_{j,i} \,t_{j,i},
\end{equation}
where the corresponding number
\begin{equation} \label{eq-Mpr}
M' = \sum_{i=0}^{k}  \epsilon_{j,i} \,v_{j,i}
\end{equation}
is used when expressing $M$. We rewrite Equations (\ref{eq-Npr}) and
(\ref{eq-Mpr}) as 
\begin{equation} \label{eq-Nprsumtj}
N' = t_{j_1} + t_{j_2} + t_{j_3} + \cdots , ~~\mbox{ and }
\end{equation}
\begin{equation} \label{eq-Mprsumvj}
M' = v_{j_1} + v_{j_2} + v_{j_3} + \cdots , ~~~~~~~~~
\end{equation}
where $j_1 > j_2 > j_3 > \cdots$. In other words, we include only the
terms for which $\epsilon_{j,i}=1$.

From Equation (\ref{eq-Nprsumtj}), we get $t_{j_1} \, \leq \,
N'$. Hence the greedy algorithm will use $t_{j_1}$ unless the next
larger number in the sequence $t_{j_1+1} \, \leq \, N'$, in which case
$t_{j_1+1}$ would be used. By the definition, $t_{j_1+1} = t_{j_1} \,
+ \, s_q$ for some $s_q$. Also, $v_{j_1+1} \, = \, v_{j_1} \, + \,
u_q$, but $v_{j_1+1}$ was not used in expressing $M'$. Hence it must
be true that $v_{j_1+1} \, = \, v_{j_1} \, + \, u_q \, > \, M'$,
whereas $t_{j_1} \, + \, s_q \, \leq \, N'$.

Now we replace the $t_j$ and $v_j$ by their corresponding combinations
of the $s_j$ and $u_j$, respectively, as follows.
\begin{equation} \label{eq-Nprminustj1}
~s_q \, \leq \, N' - t_{j_1} ~ = \, t_{j_2} + t_{j_3} + \cdots~ = ~ b_1
s_{i_1} + b_2 s_{i_2} + \cdots, ~\mbox{ and }
\end{equation}
\begin{equation} \label{eq-Mprminusvj1}
u_q \, > \, M' - v_{j_1}  = \, v_{j_2} + v_{j_3} + \cdots ~= ~
      b_1 u_{i_1} + b_2 u_{i_2} + \cdots ~ . ~~~~~~~~
\end{equation}
Since $q$ is larger than any of the $i_r$ in Equations
(\ref{eq-Nprminustj1}) and (\ref{eq-Mprminusvj1}), and since the
$\,s_i\,$ were chosen so that $\, s_i/s_{i+1} \, < \, u_i/u_{i+1}$,
from Equations (\ref{eq-Nprminustj1}) and (\ref{eq-Mprminusvj1}) we
have
\[
1 \, \leq \, b_1 (s_{i_1}/s_q) + b_2 (s_{i_2}/s_q) +  \cdots  \, < \,
    b_1 (u_{i_1}/u_q) + b_2 (u_{i_2}/u_q) +  \cdots \, < \, 1,
\]
which is a contradiction. Hence the greedy algorithm will indeed use
$t_{j_1}$.

\section{Disguising Methods} \label{sec-disgmtds}

As described above, the plaintext message $M$ can be expressed either
as $\sum d_i u_i$ or as $\sum \epsilon_{j,i} v_{j,i}$ where
$\epsilon_{j,i} = 0$ or $1$, and can be calculated using a greedy
algorithm. The sequence $\{s_i\}$ is chosen with the related auxiliary
sequence $\{ t_{j,i} \}$, which corresponds to the $v_{j,i}$. Then, if
$M$ is expressed as in Equation (\ref{eq-Mrepruivij}) and $N$ as in
Equation (\ref{eq-N}), the creator of the code can calculate the
$\epsilon_{j,i}$ from knowing $N$, and can thus compute $M$. A
disguising method that maps $s_i$ into $w_i$ is invertible if given
$T$ as expressed in Equation (\ref{eq-T}), she can compute the number
$N$. The $y_{j,i}$ are defined in terms of the $w_i$, and in the same
way $t_{j,i}$ are defined in terms of $s_i$ and $v_{j,i}$ in terms of
$u_i$.

Let $E = \max \,\sum \epsilon_{j,i}$, where the maximum is taken over
all expressions of possible messages $M$. If $n$ is the number of the
$u_i$ then $E \leq n$ since using a greedy algorithm as described, at
most one $v_{j,i}$ in each of $n$ groups is used. With more careful
analysis we can show that $E$ is much smaller than $n$, but that
result will not be critical for our analysis.

For example, in the usual modular multiplication, $w_i \equiv c
s_i\imod{m}$ or $s_i \equiv \bar{c} w_i\imod{m}$. Then
\[
\bar{c} T \equiv \sum_i d_i \bar{c} w_i \equiv \sum_i d_i s_i \equiv
N \imod {m}.
\]

The number $N$ is uniquely determined if $N < m$, which is true if $m
> E \max \{t_{j,i}\}$. Thus, it suffices to take $m > E s_{n}$. Since
the $w_i$ are defined modulo $m$, the $w_i$ are larger than the $s_i$
by a factor of $E s_n$. As such, the $w_i$ require $\log_2 (E s_n)$
bits to express in base $2$. This expansion in the size of the
disguised weights by $\log_2 (E)$ bits turns out to be similar in each
stage of the disguising.

As described earlier in Section \ref{ssec-tradknap}, the nature of
modular multiplication and the fact that some of the $s_i$ are very
small give the cryptanalyst a way to possibly compute the parameters
$c$ and $m$. We now describe some alternative disguising methods and
indicate why they could not be reversed by the cryptanalyst.

First, perform an ordinary modular multiplication so that all the
obtained values are of roughly the same size. Take any collection of
pairwise prime moduli. For simplicity, we could use distinct primes
$p_k$ for $1 \leq k \leq r$. Replace the weight $w_i$ by the
$r$-vector $\vW_i = \begin{bmatrix} w_i \imod{p_1} & w_i \imod{p_2} &
  \cdots & w_i \imod{p_r} \end{bmatrix}^T$. Then the $j^{\mbox{th}}$
component of
\begin{equation} \label{eq-jdiWi}
\sum_i d_i \vW_i \equiv \sum_i d_i w_i \imod{p_j}.
\end{equation}
Hence we can use the Chinese Remainder Theorem to compute $\sum_i d_i
w_i \imod{p_1 \dots p_k}$, and so the numerical value of $\sum_i d_i
w_i$ is determined if $p_1 \dots p_k > AE \max\{w_i\}$. Again, the
number of bits needed to express the vector $\vW_i$ is $\log_2 ( A E
\max\{w_i\})$ more than the number needed to express the $w_i$.

We can do two stages of this type of mapping, the first one with large
moduli $m_1$ and $m_2$ resulting in vectors of the form $[w'_i~
  w''_i]$ for the original weight $w_i$. The second stage uses a large
number of small primes separately on $w'_i$ and $w''_i$. These primes
could be the same or different for $w'_i$ and $w''_i$. In fact, we
could simply take $p_1=2, p_2=3, p_3=5$, and so on.

Since each of the two residues modulo $m_1$ and $m_2$ are disguised in
this way, we have two lists of residues modulo $2$ and $3$ and $5$,
and so no. The creator of the code can choose a secret permutation of
these residues. When there are $k$ primes $p_1, \dots, p_k$, each
weight is a vector of dimension $2k$. A cryptanalyst could easily see
from the size of the components whether a given component is modulo
$2$ or modulo $3$, and so on. But he cannot know which of the two
residues modulo $p_i$ came from the $m_1$ branch and which from the
$m_2$ branch of the disguising method. There are $2^k$ possible
choices here. Although it is easy to use the Chinese Remainder Theorem
on any such choice, making even one incorrect choice will produce
incorrect residues modulo $m_1$ and $m_2$. Further, the values $m_1$
and $m_2$ are not known to the cryptanalyst. Thus the security of the
disguising rests on the large number of permutations of the residues
modulo $p_i$, and not on the difficulty of solving a particular type
of calculation.

This method can be combined with other mappings as well. Another
strong candidate is using modular multiplications in the rings of
algebraic integers \cite{MoPhD2010}. This step also turns each
ordinary integer weight into a vector of dimension $k$, where $k$ is
the degree of the algebraic extension used. A sequence of such steps
can be represented as a tree. Each sequence of steps corresponds to a
different tree, a different set of parameters, and a different system
of equations describing these parameters. Even if we assume that a
cryptanalyst could obtain information about the disguising mapping
from a system of such equations, he does not even know what system of
equations is the correct one to solve. The number of such
possibilities could be made as large as we wish. Note that the final
public weight vector consists of the leaves of the corresponding tree,
which could also be secretly permuted. Thus the cryptanalyst does not
know the correct permutation nor the correct tree model to use. Any
incorrect choice results in a nondecoding result, as a consequence of
the Chinese Remainder Theorem.

\section{Cryptanalysis} \label{sec-cryptanal}

The disguising methods we just described makes an attack that hopes to
reverse the disguising steps quite unlikely to succeed. But even for
the traditional knapsack code, attacks that try to solve Equation
(\ref{eq-encodedM}) directly using basis reduction algorithms
\cite{LO85,CJLOSS92} have posed the greater threat to security.

Solving Equation (\ref{eq-encodedM}) along with the constraints that
$\epsilon_i = 0$ or $1$ could be modeled exactly using integer
optimization methods. While directly solving these instances with just
this single constraint could be difficult even for moderate values of
$n$, basis reduction-based reformulations could be more effective
\cite{KrPa2009}. But adding more constraints makes the integer
optimization instances increasingly hard to solve even when we have
only a few hundred weights, e.g., see the recent work on basis
reduction-based methods to solve market split problems
\cite{Vo2012}. Solving Equation (\ref{eq-T}) with very complex
constraints and thousands of weights appears impossible using current
methods.

Direct basis reduction-based approaches, as opposed to integer
optimization approaches, can handle much larger problems. But these
methods cannot impose strict adherence to the constraints as the
integer optimization models do. Instead, these algorithms find short
vectors in appropriately defined lattices, which correspond to
solutions of Equations like (\ref{eq-encodedM}) and (\ref{eq-T}). In
the case of Equation (\ref{eq-encodedM}), the shortest vector in the
lattice is the desired vector of $\epsilon_i$, even though the basis
reduction-based methods are not guaranteed to find this particular
vector. Indeed, the {\em shortest vector problem} (SVP) and the
closely related {\em closest vector problem} (CVP) are known to be
hard problems. CVP is known to be NP-complete, and so is SVP under
randomized reductions \cite{MiGo2002}. Still, such algorithms are
often successful in practice to solve the problem instances exactly
\cite{GamNgupredlatred08}.

Most basis reduction-based approaches on the default knapsack code in
Equation (\ref{eq-encodedM}) start by defining an appropriate lattice
in which the shortest vector corresponds to the correct decoding
message. For the subset sum problem with weights $w_i$ and the target
sum $T$, Coster et al.~\cite{CJLOSS92} consider the lattice generated
by
\begin{equation} \label{eq-latt}
  \mathcal{L} = \begin{bmatrix}
    c \, \vw^T & c \, T \\
    2I   & \vone \end{bmatrix},
\end{equation}
where $\vw$ is the vector of weights $w_i$, $I$ is the $n \times n$
identity matrix, and $\vone$ is the $n$-vector of ones. When the
multiplier $c$ is chosen large enough, the vector $\vep$ corresponding
to the correct decoding will generate the shortest vector in this
lattice by multiplying the first $n$ columns of $\mathcal{L}$ by
$\vep$ and the last column by $-1$. This shortest vector has a length
of $\sqrt{n}$. To locate this shortest vector, one tries to find
short(est) vectors in the lattice $\mathcal{L}$. While it is not
guaranteed to find the correct decoding vector in every case, these
algorithms succeed with high probability when the {\em density} of the
knapsack, defined as $\, n/(\max_i \log(w_i))$, is small. When $n$ is
of the order of a few hundreds, these methods have been shown to be
effective in finding the correct decoding \cite{HJ10}.

To see how these methods might be applied to our code, we examine the
example described in Section \ref{ssec-exmplcode}. The cryptanalyst
has the choice of trying to solve either
\begin{equation} \tag{(\ref{eq-T}) revisited}
T = \sum_{i=0}^{n-1} d_i w_i \mbox{\hspace*{0.5in}or \hspace*{0.5in}}
T = \sum_{i=0}^{n-1} \epsilon_{j,i} \, y_{j,i}.
\end{equation}
We examine both possibilities by first calculating the expected
lengths of the vectors $\vd = [d_0~\cdots~d_{n-1}]^T$ of dimension
$n=1000$ and the vector $\vep$ of $\epsilon_{i,j}$ values of length
$10n = 10^4$. It is not clear how one would add constraints forbidding
the particular substrings of digits that are not permitted in $\vd$ in
the former approach. While one could potentially model all constraints
forbidding nonallowed combinations of $\epsilon_{i,j}$ in the latter
approach, this step would produce a candidate lattice $\mathcal{L}$ in
which the single knapsack row is replaced by a substantially large
number of simultaneous linear Diophantine equations. Hence the
original basis reduction-based methods will struggle to find the
correct decoding vector in this case as well. We could add one further
level of difficulty in modeling the correct lattice $\mathcal{L}$ for
the latter approach using $\epsilon_{i,j}$. We could slightly alter
$M$ by adding a small number $\bar{M}$, or by using $2M$ or $3M$
instead of $M$, so that the expected length of the decoding vector
becomes longer. The value $\bar{M}$, or the multiplier, is sent in the
clear, and should not affect the security of the code.

Further, since the first few groups of the $v_{j,i}$ do not have the
same number of elements, it is convenient to start with, say, $u_{20}$
instead of $u_0$. In the encoding process, the remaining values
smaller than $u_{20}$ are also sent in the clear. This modification
also makes the choice of the sequence of $s_i$ easier to make, without
increasing the size of the final public weights much. Sending these
small values in the clear does not change the efficiency of the code
significantly.

In the encoding procedure, suppose at some stage we have the value
$M'$ as in Equation (\ref{eq-Mpr}). If $M'$ falls randomly in the
interval $[u_k, u_{k+1})$, we first calculate the probability that
  $v_{j,k} \leq M' < v_{j+1,k}$. Then we calculate the probability
  that $M' - v_{j,k}$ falls in a particular smaller interval
  $[u_{\ell}, u_{\ell+1})$. Recall that the associated principal
    eigenvalue is $\alpha \approx 1.9754$, and we may approximate
    $u_k$ by $\alpha^k$ for large $k$. There are $10$ subintervals to
    consider. Letting $v_{j,k} = v_j$ temporarily for ease of
    notation, the intervals are $I_1 = [v_1, \, v_2 = v_1 + u_{k-2})$,
      $I_2 = [v_2, \, v_2 + u_{k-3})$, $I_3 = [v_3, \, v_3 +
          u_{k-3})$, and $I_j = [v_j, \, v_j + u_{k-4})$ for $4 \leq j
            \leq 10$, where $v_{10} + u_{k-4} = u_{k+1}$.

If $h \leq k-5$, then $M' - v_{j,k} \, \in \, [u_h, u_{h+1})$ can
  occur in any of the $10$ subintervals. If $h=k-4$, then $M' -
  v_{j,k} \in [u_{k-4}, u_{k-3})$ can occur only in intervals $I_1,
    I_2$, or $I_3$, and $M' - v_{j,k} \in [u_{k-3}, u_{k-2})$ can
      occur only in $I_1$. Hence the probability that $M' - u_k \in
      [u_{k-5}, u_{k-4})\,$ is $\,10 (u_{k-4} - u_{k-5})\,/\,(u_{k+1}
        - u_k) \approx 10 \alpha^{-5}$, while the probability that $M'
        - u_k \in [u_{k-4}, u_{k-3})\,$ is $\,3 (u_{k-3} -
          u_{k-4})\,/\,(u_{k+1} - u_k) \approx 3 \alpha^{-4}$.
If $M' \in [u_k, u_{k+1})$ and $M' - u_k \,\in\, [u_h, u_{h+1})$, then
there are $k-h-1$ groups of the $v_{j,i}$ that are skipped, i.e., they
are not present in the representation of $M'$. Then the expected
number of skipped groups between each pair that does appear is
approximately
\begin{align*}
        & ~10 \sum_{h=0}^{k-5} \alpha^{h-k} \,(k-h-1) + 3 \alpha^{-4} (3) +  \alpha^{-3}(2)  \\
  =     & ~10 \sum_{j=5}^k \alpha^{-j} \, (j-1) + 9 \alpha^{-4} + 2 \alpha^{-3} \\
\approx & ~10 \left(\alpha^{-3} + 3 (\alpha-1) \alpha^{-4} \,/\,(\alpha-1)^2 \right) + 9 \alpha^{-4} + 2 \alpha^{-3} \\ 
\approx & ~~~3.383 ~+~ 0.591 ~+~ 0.259 ~~~= ~4.233.
\end{align*}
Therefore the expected number of groups represented when expressing
$M$ is $n\,/\,5.233$, or $191$ when $n=1000$.

In this example, $M \approx \alpha^{1000} \approx 5.5 \times 10^{295}
\approx 2^{982}$. After the disguising steps described above, the size
of the target sum $T$ in number of bits needed is about $1030$. There
are then $2^{1030} \approx 10^{310}$ possible target objects or
vectors. Even if we use extra disguising steps, the number of target
objects might be as large as $10^{325}$.

We expect the representation of $M$ to use some $v_{i,j}$ from $191$
of the $1000$ groups. As described above, we can easily alter $M$
slightly to make sure at least $191$ groups are used, i.e., the vector
corresponding to the correct decoding has length at least $191$. If we
look at shorter vectors, say of length $180$, consisting of one
$v_{i,j}$ from $180$ different groups, there are
\[
\nchoosek{1000}{180} \, 10^{180} ~ \approx ~ 10^{383}
\]
such vectors. Thus we expect $10^{383-325} = 10^{58}$ of these vectors
to yield the same $T$. There are many even shorter vectors that yield
$T$. If a basis reduction-based method finds the shortest vector
yielding $T$, it will not be the one needed to decode $T$. If the
method finds some short vector at random, the chance that it is the
correct decoding vector would be much smaller than $10^{-58}$.

\medskip
Now suppose the cryptanalyst tries to solve Equation (\ref{eq-T})
directly for the vector of digits $\vd$ instead. The candidate vectors
have dimension $1000$ now, but the entries may be as large as $6$. 

To compute the expected length of $\vd$, we need to calculate the
expected value of each digit. By the calculations described above, we
expect $191$ of the $10$ blocks $100, 1010, \dots, 10126$ to
appear. The $j^{\mbox{th}}$ block is used if $M$ is in the interval
$I_j$ described above. Hence the block $100$ is used with probability
$\, u_{k-2} \, / \,(u_{k+1} - u_k) \, \approx \, \alpha^{-2}
/(\alpha-1)$, $1010$ and $1011$ appear with probability
$\alpha^{-3}/(\alpha-1)$, and the other $7$ blocks with probability
$\alpha^{-4}/(\alpha-1)$.

Among the $1000$ digits $d_i$, we expect roughly $475$ zeros, $370$
ones, $103$ twos, and $13$ each of $3,4,5$, and $6$. The expected
length of the vector is $1900$. There are more than $10^{441}$ vectors
of length $1550$ consisting of $200$ zeros, $500$ ones, $150$ twos and
$50$ threes. There exist even more shorter vectors. Since there are
only $10^{325}$ target vectors, we expect many shorter vectors to
correspond to the same value $T$. Thus, unless it is possible to force
the basis reduction-based algorithms to exclude vectors that do not
correspond to allowable strings of digits, these methods will yield
shorter vectors that do not correctly decode the encrypted message
$T$. Integer optimization-based methods could exclude nonallowed
vectors, but they could not handle problems of this size.

The example code presented above is actually smaller and simpler than
one that would be suggested for use. Indeed, one could take $n$ as
large as $10^5$ without needing too large an amount of memory, and $A$
as large as $100$ or $1000$. This setting would require the
cryptanalyst to solve problems with $10^8$ vectors if he used the
version of Equation (\ref{eq-T}) with $y_{j,i}$.

\section{Discussion} \label{sec-disc}

We have described how to use recurrence sequence representations to
create a knapsack-like public key code. There is considerable freedom
in the choice for the encoding sequence $\{u_i\}$, as well as its
order $h$, the size and pattern of the coefficients $a_j$, and the
number of terms $n$. Larger orders and coefficients create more false
short vectors and hence increase the complexity of
cryptanalysis. However, they also increase the size of $u_i$ and the
weights $w_i$. Questions about optimal choices for these parameters
remain for further study. 

Another relevant question is whether specific number theoretic
properties of such representations could be used to make this code
even more secure. In particular, security against attacks targeting
the representation using $\epsilon_{i,j}$ relies of the size of the
instances, which could not be handled by state of the art methods. But
could we make such attacks less feasible, or even impossible to mount,
using inherent properties of the sequence used for the representation?

\bibliographystyle{plain} 

\bibliography{IP_Crypt_Refs}

\end{document}